\documentclass[12pt]{amsart}
\usepackage{amssymb}
\usepackage{graphics}
\usepackage{latexsym}
\usepackage{amsmath}
\usepackage{amssymb,amsthm,amsfonts}
\usepackage{amscd}
\usepackage[arrow, matrix, curve]{xy}
\usepackage{syntonly}
\title[Abelian variety of Mumford's type]{On the Newton polygons of abelian varieties of Mumford's type}
\author[Mao Sheng]{Mao Sheng}
\author[Kang Zuo]{Kang Zuo}
\email{sheng@uni-mainz.de} \email{zuok@uni-mainz.de}
\address{Institut f\"{u}r  Mathematik, Universit\"{a}t
Mainz, Mainz, 55099, Germany} \iffalse
\begin{document}
%%%%%%%%%%%%%%%%%%%% Text italic %%%%%%%%%%%%%%%%%%%%%%%%%%%%
\theoremstyle{plain}
\newtheorem{thm}{Theorem}[section]
\newtheorem{theorem}[thm]{Theorem}
\newtheorem{lemma}[thm]{Lemma}
\newtheorem{corollary}[thm]{Corollary}
\newtheorem{proposition}[thm]{Proposition}
\newtheorem{addendum}[thm]{Addendum}
\newtheorem{variant}[thm]{Variant}
%%%%%%%%%%%%%%%%%%%% Text roman %%%%%%%%%%%%%%%%%%%%%%%%%%%%%
\theoremstyle{definition}
\newtheorem{lemma and definition}[thm]{Lemma and Definition}
\newtheorem{construction}[thm]{Construction}
\newtheorem{notations}[thm]{Notations}
\newtheorem{question}[thm]{Question}
\newtheorem{problem}[thm]{Problem}
\newtheorem{remark}[thm]{Remark}
\newtheorem{remarks}[thm]{Remarks}
\newtheorem{definition}[thm]{Definition}
\newtheorem{claim}[thm]{Claim}
\newtheorem{assumption}[thm]{Assumption}
\newtheorem{assumptions}[thm]{Assumptions}
\newtheorem{properties}[thm]{Properties}
\newtheorem{example}[thm]{Example}
\newtheorem{conjecture}[thm]{Conjecture}
\numberwithin{equation}{thm}

% Skriptbuchstaben
\newcommand{\pP}{{\mathfrak p}}
\newcommand{\sA}{{\mathcal A}}
\newcommand{\sB}{{\mathcal B}}
\newcommand{\sC}{{\mathcal C}}
\newcommand{\sD}{{\mathcal D}}
\newcommand{\sE}{{\mathcal E}}
\newcommand{\sF}{{\mathcal F}}
\newcommand{\sG}{{\mathcal G}}
\newcommand{\sH}{{\mathcal H}}
\newcommand{\sI}{{\mathcal I}}
\newcommand{\sJ}{{\mathcal J}}
\newcommand{\sK}{{\mathcal K}}
\newcommand{\sL}{{\mathcal L}}
\newcommand{\sM}{{\mathcal M}}
\newcommand{\sN}{{\mathcal N}}
\newcommand{\sO}{{\mathcal O}}
\newcommand{\sP}{{\mathcal P}}
\newcommand{\sQ}{{\mathcal Q}}
\newcommand{\sR}{{\mathcal R}}
\newcommand{\sS}{{\mathcal S}}
\newcommand{\sT}{{\mathcal T}}
\newcommand{\sU}{{\mathcal U}}
\newcommand{\sV}{{\mathcal V}}
\newcommand{\sW}{{\mathcal W}}
\newcommand{\sX}{{\mathcal X}}
\newcommand{\sY}{{\mathcal Y}}
\newcommand{\sZ}{{\mathcal Z}}
% Sonderbuchstaben mit Doppellinie
\newcommand{\A}{{\mathbb A}}
\newcommand{\B}{{\mathbb B}}
\newcommand{\C}{{\mathbb C}}
\newcommand{\D}{{\mathbb D}}
\newcommand{\E}{{\mathbb E}}
\newcommand{\F}{{\mathbb F}}
\newcommand{\G}{{\mathbb G}}
\newcommand{\HH}{{\mathbb H}}
\newcommand{\I}{{\mathbb I}}
\newcommand{\J}{{\mathbb J}}
\renewcommand{\L}{{\mathbb L}}
\newcommand{\M}{{\mathbb M}}
\newcommand{\N}{{\mathbb N}}
\renewcommand{\P}{{\mathbb P}}
\newcommand{\Q}{{\mathbb Q}}
\newcommand{\R}{{\mathbb R}}
\newcommand{\SSS}{{\mathbb S}}
\newcommand{\T}{{\mathbb T}}
\newcommand{\U}{{\mathbb U}}
\newcommand{\V}{{\mathbb V}}
\newcommand{\W}{{\mathbb W}}
\newcommand{\X}{{\mathbb X}}
\newcommand{\Y}{{\mathbb Y}}
\newcommand{\Z}{{\mathbb Z}}
\newcommand{\id}{{\rm id}}
\newcommand{\rank}{{\rm rank}}
\newcommand{\END}{{\mathbb E}{\rm nd}}
\newcommand{\End}{{\rm End}}
\newcommand{\Hom}{{\rm Hom}}
\newcommand{\Hg}{{\rm Hg}}
\newcommand{\tr}{{\rm tr}}
\newcommand{\Cor}{{\rm Cor}}
\newcommand{\Res}{{\rm Res}}
\newcommand{\GL}{\mathrm{GL}}
\newcommand{\GSp}{\mathrm{GSp}}
\newcommand{\Sp}{\mathrm{Sp}}
\newcommand{\SL}{\mathrm{SL}}
\newcommand{\Aut}{\mathrm{Aut}}
\newcommand{\Sym}{\mathrm{Sym}}
\newcommand{\Gal}{\mathrm{Gal}}
\newcommand{\DD}{\mathbf{D}}
\newcommand{\SU}{\mathrm{SU}}
%%%%%%%%%%%%%%%%%%%%%%%%%%%%%%%%%%%%%%%%%%%%%%%%%%%%%%%%
\maketitle

\footnotetext[1]{This work was supported by SFB/Transregio 45
Periods, Moduli Spaces and Arithmetic of Algebraic Varieties of the
DFG (Deutsche Forschungsgemeinschaft).}

\begin{abstract}
Let $A$ be an abelian variety of Mumford's type. This paper
determines all possible Newton polygons of $A$ in char $p$. This
work generalizes a result of R. Noot in \cite{Noot1}-\cite{Noot2}.
\end{abstract}

\section{Introduction}
In \cite{Mu} D. Mumford constructs families of abelian varieties
over smooth projective arithmetic quotients of the upper half plane,
whose general fiber has only $\Z$ as its endomorphism ring. This is
the first example of Shimura families of Hodge type. Roughly
speaking, a Shimura family of Hodge type is characterized by the
Hodge classes on the tensor products of the weight one
$\Q$-polarized Hodge structure of an abelian variety and its dual.
When they are generated by $(1,1)$-classes, it is called a Shimura
curve of PEL type: By the Lefschetz $(1,1)$ theorem, they are
classes of algebraic cycles. Since the Hodge conjecture or its alike
has not been settled in general, it is often the case that a
result can be established relatively easily for the PEL type than for the Hodge type.\\

In a series of papers \cite{Noot0},\cite{Noot1},\cite{Noot2} R. Noot
studies the Galois representations associated to an abelian variety
over a number field which appears as a closed fiber of a Mumford's
family, and he obtains various results about such an abelian
variety, notably potentially good reduction, classification and
existence of isogeny types and Newton polygons.
Some of the arguments in loc. cit. work also in a general setting.\\

From a completely different direction, E. Viehweg and the second
named author characterize those families of abelian varieties which
reach the Arakelov bound (cf. \cite{VZ}), and it turns out that
under a natural assumption one can classify them up to isogeny and
finite \'{e}tale base change. The results are closely related with
Mumford's families and we restate them here as follows:
\begin{theorem}[Viehweg-Zuo, Theorem 0.2, Theorem 0.5,
\cite{VZ}]\label{theorem of VZ} Let $f: X\to Y$ be a semi-stable
family of $g$-dimensional abelian varieties, smooth over $U:=Y-S$.
Assume that there is no unitary local subsystem in
$R^1f_{*}\C_{f^{-1}U}$. Then if $f$ reaches the Arakelov bound,
namely the following Arakelov inequality
$$
2.\deg(f_{*}\Omega_{X|Y}(\log f^{-1}(S)))\leq g\cdot
\deg(\Omega^1_{Y}(\log S))
$$
becomes equality, then after a possible finite \'{e}tale base
change, $f$ is isogenous to
\begin{itemize}
    \item [{\bf Case 1:}] the $g$-fold self product of a universal semi-stable family
    of elliptic curves over a modular curve if $S\neq \emptyset$,
    \item [{\bf Case 2:}] the self product of a Shimura family of
    Mumford type if $S= \emptyset$.
\end{itemize}
\end{theorem}
The latter case generalizes the example of Mumford in a natural way.
The construction in loc. cit. will be briefly recalled (see \S1 for
more details): Let $F$ be a totally real field of degree $d\geq 1$
and $D$ a quaternion division algebra over $F$, which is split only
at one real place of $F$. The corestriction $\Cor_{F|\Q}D$ is
isomorphic to $M_{2^d}(\Q(\sqrt{b}))$ for a rational number $b$.
Depending on whether $b$ is a square or not, one is able to
construct Shimura families of Hodge type with the fiber dimension
$2^{d+\epsilon(D)-1}$. The cases $d=1,2$ give rise to Shimura
families of PEL type, and the case $d=3$ is just the
example of Mumford.\\

By the theory of canonical models of Shimura varieties, these
families are defined naturally over certain number fields. By abuse
of notation, we call a $\bar \Q$-closed fiber of such a Shimura
family an abelian variety of \emph{Mumford's type}. The purpose of
the paper is to study the $p$-adic Galois representation of an
abelian variety of Mumford's type. That is, we intend to study the
$p$-adic Hodge structure of such an abelian variety using the
$p$-adic Hodge theory of Fontaine and his school. Certainly there
will be some overlaps of our approach with what Noot has done in his
study on Mumford's example (see particularly \S3 \cite{Noot1}, \S2
\cite{Noot2}). However we shall emphasize that the relation between
two dimensional potentially crystalline $\Q_{p^r}$-representations
and $p$-adic Hodge structures of abelian varieties of Mumford's type
is essentially new. A first study of this relation gives us the
classification of Newton polygons of good reductions of abelian
varieties of Mumford's type. To make the result complete, we also
include the existence of possible Newton polygons, which is based on
Noot \S3-5 \cite{Noot2}.

\begin{theorem}\label{classification and existence of newton polygon}
Let $F$ be a totally real field of degree $d\geq 1$ and $D$ a
quaternion division algebra over $F$, which is split only at one
real place of $F$. Let $A$ be an abelian variety of Mumford's type
associated to $D$, which is defined over a number field $K\supset
F$. Let $p$ be a rational prime number satisfying Assumption
\ref{assumption on p} and $\mathfrak{p}_K$ a prime of $K$ over $p$.
Assume that $A$ has the good reduction $A_k$ at $\mathfrak{p}_K$.
Put $\mathfrak{p}=\mathfrak{p}_K\cap \sO_F$ and
$r=[F_{\mathfrak{p}}:\Q_p]$. Then the Newton polygon of $A_k$ is
either $\{2^{d+\epsilon(D)}\times \frac{1}{2}\}$ (i.e.
supersingular) or
$$
\{2^{d-r+\epsilon(D)}\times 0, \cdots, 2^{d-r+\epsilon(D)}\cdot{r
\choose i}\times \frac{i}{r}, \cdots, 2^{d-r+\epsilon(D)}\times 1\}.
$$
Here $\epsilon(D)$ is equal to 0 or 1 depending only on $D$ (see \S1
for the definition). Furthermore each possible Newton polygon does
occur for an abelian variety of Mumford's type associated to $D$.
\end{theorem}

The above result confirms the conjecture on the Newton
stratification of a Shimura curve of
Hodge type made in \cite{SZZ} partially.\\

{\bf Acknowledgements:} The paper is indebted to the initial works
of Rutger Noot \cite{Noot1}-\cite{Noot2}. It is his formula
Proposition 2.2 \cite{Noot2} that inspires the current work. We
thank him heartily. We would like to thank Jean-Marc Fontaine for
several helpful discussions on this paper. He pointed out to us that
the tensor factor $V_{1}$ can be crystalline and even wrote down one
of the possible Newton polygons, as what we have managed to prove in
Theorem \ref{possible newton slopes}. We owe the existence of the
paper to this insight. We would like also to express our sincere
thanks to Yi Ouyang and Liang Xiao for their helps on the proof of
Theorem \ref{HT implies de-rham }, that should be due to Laurent
Berger or his student Giovanni Di Matteo.

\section{Abelian varieties of Mumford's type and tensor decomposition of Galois representations}
Let us start with the following
\begin{lemma}[Lemma 5.7 (a) \cite{VZ}]
Let $F$ be a totally real field of degree $d$ and $D$ be a
quaternion division algebra over $F$, which is split at a unique
real place of $F$. Let $\Cor_{F|\Q}(D)$ be the corestriction of $D$
to $\Q$. Then
\begin{itemize}
    \item [(i)] $\Cor_{F|\Q}(D)\cong M_{2^d}(\Q)$ and $d$ is odd,
    or
    \item [(ii)] $\Cor_{F|\Q}(D)\ncong M_{2^d}(\Q)$. Then
    $$
\Cor_{F|\Q}(D)\cong M_{2^d}(\Q(\sqrt{b})),
    $$
    where $\Q(\sqrt{b})$ is a(n) real (resp. imaginary) quadratic
    field extension of $\Q$ if $d$ is odd (resp. even).
\end{itemize}
\end{lemma}
For a given $D$, an isomorphism in the above lemma will be fixed
once for all. Clearly we can write either case $\Cor_{F|\Q}(D)\cong
M_{2^d}(\Q(\sqrt{b}))$ for a rational number $b\in \Q$. If we
require $b$ to be square free, then $b$ is uniquely determined by
$D$.  We define a function $\epsilon(D)$ which takes the value 0
when $b=1$ and 1 otherwise. For case (ii) we shall also further
choose an embedding $M_{2^d}(\Q(\sqrt{b}))\hookrightarrow
M_{2^{d+1}}(\Q)$. So we shall fix an embedding
$\Cor_{F|\Q}(D)\hookrightarrow
M_{2^{d+\epsilon(D)}}(\Q)$ in either case.\\

Let $\bar{\ }$ be the standard involution of $D$. One defines a
$\Q$-simple group $$\tilde{G'}:=\{x\in D|\ x\bar x=1\},$$ and
$\tilde G:=\G_{m,\Q}\times \tilde{G'}$. Recall that there is a
natural morphism of $\Q$-groups (cf. \S4 \cite{Mu})
$$
{\rm Nm}: D^*\to \Cor_{F|\Q}(D)^*\hookrightarrow
\GL_{2^{d+\epsilon(D)}}(\Q).$$ One defines a reductive $\Q$-group
$G$ to be the image of the morphism $ \tilde G\to
\GL_{2^{d+\epsilon(D)}}(\Q)$, that is the product of the natural
morphism $\G_{m,\Q}\to \GL_{2^{d+\epsilon(D)}}(\Q)$ and ${\rm
Nm}|_{\tilde{G'}}$. The image of $\tilde{G'}$ in $G$ is denoted by
$G'$. The resulting natural morphism
$$N: \tilde G=\G_{m,\Q}\times
\tilde{G'}\to G$$ is a central isogeny. The set of real places of
$F$ is denoted by $\{\tau_i\}_{1\leq i\leq d}$, and we assume that
$D$ is split over $\tau:=\tau_1$. Then one has an isomorphism of
real groups:
$$
\tilde{G'}(\R)\simeq \SL_{2}(\R)\times \SU(2)^{\times d-1}.
$$
One defines
$$
u_0: S^1\to \tilde{G'}(\R),\ e^{i\theta}\mapsto \left(%
\begin{array}{cc}
  \cos \theta & \sin \theta \\
  -\sin \theta & \cos \theta \\
\end{array}%
\right)\times id^{\times d-1},
$$
and $\tilde{h}_0=id\times u_0: \R^*\times S^1\to \tilde{G}(\R)$
descends to a morphism of real groups:
$$
h_0: \SSS\to G_{\R}.
$$
Let $X$ be the $G(\R)$-conjugacy class of $h_0$, and one verifies
that $(G,X)$ gives the Shimura datum of a Shimura curve. Let
$W_{\Q}=\Q(\sqrt{b})^{2^d}$ with the natural action of $G'$. It is
easy to verify that there exists a unique $\Q(\sqrt{b})$-valued
symplectic form $\omega$ on $W_{\Q}$ up to scalar invariant under
$G'$-action. One defines a $\Q$-valued symplectic form
$\psi=Tr_{\Q(\sqrt{b})|\Q}\omega$ on $W_{\Q}$ and
$(\GSp(W_{\Q},\psi),X(\psi))$ be the Siegel modular variety. Then
there is a natural embedding $G\hookrightarrow \GSp(W_{\Q},\psi)$
which carries $X$ into $X(\psi)$: This realizes $(G,X)$ as a Shimura
curve of Hodge type. Let $C\subset G(\A_f)$ be a compact open
subgroup and one defines the Shimura curve as the double coset
$$
Sh_{C}(G,X):= G(\Q)\backslash X\times G(\A_f)/ C,
$$
where
$$
q(x,a)b=(qx,qab), \ q\in G(\Q), x\in X, a\in G(\A_f), b\in C.
$$
The theory of canonical model says that $Sh_{C}(G,X)$ is naturally
defined over its reflex field $\tau(F)\subset \C$. By choosing $C$
small enough, there exists a universal family of abelian varieties
$\sX_{C}\to Sh_{C}(G,X)$ which is also defined over $\tau(F)$. A
Shimura family of Mumford type in Theorem \ref{theorem of VZ} is a
geometrically connected component of such a universal family
(different components are isomorphic to each other over $\bar \Q$).
We call a $\bar \Q$-closed fiber of
such a Shimura family an abelian variety of Mumford's type. \\

We take the following assumption on the prime number $p$:
\begin{assumption}\label{assumption on p}
Let $p$ be a rational prime which does not divide the discriminants
of $F$ and $D$.
\end{assumption}
we have then
$$
p\sO_F=\prod_{i=1}^{n}\mathfrak{p}_i.
$$
Let $\bar \Q$ be the algebraic closure of $\Q$ in $\C$. By choosing
an embedding $\bar \Q\to \bar \Q_p$, one gets an identification
$$
\Hom_{\Q}(F,\bar \Q)=\Hom_{\Q}(F,\bar \Q_p)=
\coprod_{i=1}^{n}\Hom_{\Q_p}(F_{\mathfrak{p}_i},\bar \Q_p).$$ In the
rest of the paper such an embedding will be fixed. By a possible
rearrangement of indices, we can assume that $\tau\in
\Hom_{\Q_p}(F_{\mathfrak{p}_1},\bar \Q_p)$ under the above
identification. By the above assumption one has an isomorphism
$$
D\otimes_{F}F_{\mathfrak{p}_1}\simeq M_2(F_{\mathfrak{p}_1}).$$
Such an isomorphism will be also fixed. Put $r=[F_{\mathfrak{p}_1}:\Q_p]$.\\

Let $A$ be an abelian variety of Mumford's type defined over a
number field $K$. Let $\rho: \Gal_K\to \GL(H^1_{et}(\bar A, \Q_p))$
be the associated Galois representation. For simplicity we write
$H_{\Q_p}$ for $H^1_{et}(\bar A, \Q_p)$. The aim of this section is
to show the following tensor decomposition:
\begin{proposition}\label{tensor decomposition}
Restricted to an open subgroup $\Gal_{K'}\subset \Gal_K$,
$(\rho,H_{\Q_p})$ has a natural tensor decomposition:
$$
H_{\Q_p} \simeq (V_{\Q_p} \otimes U_{\Q_p})^{\oplus
2^{\epsilon(D)}},
$$
together with a further tensor decomposition of $V_{\Q_p}$:
$$
V_{\Q_p}\otimes_{\Q_{p}} \Q_{p^r}\simeq
V_{1}\otimes_{\Q_{p^r}}V_{1,\sigma}\otimes_{\Q_{p^r}}
\cdots\otimes_{\Q_{p^r}} V_{1,\sigma^{r-1}},
$$
where for $0\leq i\leq r-1$, $V_{1,\sigma^{i}}$ is the
$\sigma^i$-conjugate of the two dimensional
$\Q_{p^r}$-representation $V_1$.
\end{proposition}
\begin{proof}
Let $H_{\Q}=H^1_{B}(A(\C),\Q)$ be the singular cohomology, and
$G_{A}\subset \GL(H_{\Q})$ be the Mumford-Tate group of $A$
(strictly speaking it is the projection of the Mumford-Tate group of
$A$ to the first factor of $\GL(H_{\Q})\times \G_{m,\Q}$ cf. \S3,
\cite{De}). By construction of the Shimura curve $Sh_{C}(G,X)$,
under a natural identification $H_{\Q}=W_{\Q}$ the Mumford-Tate
group $G_A$ is a subgroup of $\xi:G\hookrightarrow \GL(W_{\Q})$. By
Proposition 2.9 (b) in loc. cit., there is a finite field extension
$K\subset K'\subset \bar \Q$ such that the representation
$$\rho|_{\Gal_{K'}}:\Gal_{K'}\to \GL(H_{\Q_p})=\GL(W_{\Q_p})$$
factors through
$$\xi_{\Q_p}:
G(\Q_p)\hookrightarrow \GL(W_{\Q_p}).$$ In the following we examine
the morphism $\xi_{\Q_p}$ closely. One has a short exact sequence of
$\Q$-groups:
$$
1\to G'\to G\to \G_{m}\to 1.
$$
In order to show the tensor decomposition it suffices to consider
the restriction of $\xi$ to the subgroup $G'$ and its base change to
$\Q_p$. Recall that $\Cor_{F|\Q}(D)$ is defined as the
$\Gal_{\Q}$-invariants of the $\bar \Q$-algebra
$\otimes_{i=1}^{d}D_{\tau_i}$ where
$D_{\tau_i}=D\otimes_{F}(\bar \Q,\tau_i)$. Since $D_{\tau_i}$ is split, it acts on $H_{i}=D_{\tau_i}\cdot\left(%
\begin{array}{cc}
  1 & 0 \\
  0 & 0 \\
\end{array}%
\right)$ by the left multiplication. It is isomorphic to the
standard representation of $M_{2}(\bar \Q)$ on $\bar{\Q}^{2}$. Then
as $\Cor_{F|\Q}(D)\otimes \bar \Q$-modules one has a natural tensor
decomposition
$$
W_{\Q}\otimes \bar \Q=(\otimes_{i=1}^d H_i)^{\oplus
2^{\epsilon(D)}}.
$$
The decomposition group $D_p\subset \Gal_\Q$, defined by the chosen
embedding $\bar \Q\hookrightarrow \bar \Q_p$, is isomorphic to the
local Galois group $\Gal_{\Q_p}$. It acts on the set
$\{\tau_i\}_{1\leq i\leq d}$ and by Assumption \ref{assumption on p}
it decomposes into $n$ orbits according to the prime decomposition.
WLOG we can write $\{\tau=\tau_1,\cdots,\tau_r\}$ for the orbit
containing $\tau$. Put then $V_{\bar
\Q_p}=\otimes_{i=1}^{r}(H_{i}\otimes_{\bar \Q}\bar \Q_p)$ and
$U_{\bar \Q_p}=\otimes_{i=r+1}^{d}(H_{i}\otimes_{\bar \Q}\bar
\Q_p)$. Thus one has a natural tensor decomposition
$$
W_{\Q}\otimes \bar \Q_p= (V_{\bar \Q_p}\otimes U_{\bar
\Q_p})^{\oplus 2^{\epsilon(D)}}.
$$
Taking the $\Gal_{\Q_p}$-invariants, one obtains a tensor
decomposition as $\Cor_{F|\Q}(D)\otimes \Q_p$-modules
$$
W_{\Q_p}=(V_{\Q_p}\otimes U_{\Q_p})^{\oplus 2^{\epsilon(D)}}.
$$
Since the action of $G'(\Q_p)$ on $W_{\Q_p}$ via $\xi_{\Q_p}$
factors through that of $\Cor_{F|\Q}(D)\otimes \Q_p$, the above
decomposition is also a tensor decomposition for the $\Gal_{K'}$
action. To proceed further, we consider the base change of
$$
{\rm Nm}: D^*\to \Cor_{F|\Q}(D)^*, \ d\mapsto (d\otimes 1)\otimes
\cdots\otimes (d\otimes 1).
$$
to $\Q_p$. Let $F_i$ be the completion of $F$ at the prime
$\mathfrak{p}_i$. Then one has a natural isomorphism ${\rm
Nm}\otimes \Q_p=\prod_{i=1}^{n}{\rm Nm}_i$ with
$$
{\rm Nm}_i: (D\otimes_FF_i)^*\to \Cor_{F_i|\Q_p}(D\otimes_FF_i)^*.
$$
By Assumption \ref{assumption on p}, $F_1\simeq \Q_{p^r}$ and
$D\otimes_FF_1\simeq M_2(\Q_{p^r})$. It is clear that the action of
$\Gal_{\Q_p}$ on
$$
\Cor_{F_1|\Q_p}(D\otimes_FF_1)\otimes_{\Q_p} \bar \Q_p\simeq
[M_{2}(\Q_{p^r})\otimes _{\Q_{p^r}}(\bar
\Q_p,\tau_1)]\otimes\cdots\otimes [M_{2}(\Q_{p^r})\otimes
_{\Q_{p^r}}(\bar \Q_p,\tau_r)]
$$
factors through $$\Gal_{\Q_p}\twoheadrightarrow
\Gal_{\Q_{p^r}|\Q_p}=\langle\sigma\rangle.$$ Thus it descends to the
tensor decomposition
$$
\Cor_{F_1|\Q_p}(D\otimes_FF_1)\otimes_{\Q_p}
F_1=\otimes_{i=0}^{r-1}[M_{2}(F_1)\otimes_{F_1}(F_1,\sigma^{i})]=\otimes_{i=0}^{r-1}M_{2}(F_1^{\sigma^i})
.$$ This gives us the further tensor decomposition of $V_{\Q_p}$.
\end{proof}

\section{Two dimensional potentially crystalline $\Q_{p^r}$-representations}
A large portion of this section is expository and it is based on
\S7.3 \cite{FO}. First we recall several standard notions in the
theory of $p$-adic representations (see for example \cite{FO}). Let
$E$ be a finite field extension of $\Q_p$ and $V$ a $p$-adic
representation of $\Gal_E$. Let $B$ be one of the period rings
$B_{HT},B_{dR},B_{crys},B_{st}$ introduced by Fontaine.
\begin{definition} A representation $V$ is called
$B$-admissible if the tensor representation $V\otimes B$ of $\Gal_E$
is trivial, namely there exists a $B$-basis of $V\otimes B$ on which
$\Gal_E$ acts trivially.
\end{definition}
Let $D_B(V)=(V\otimes B)^{\Gal_E}$ which is naturally a
$B^{\Gal_E}$-vector space. A nice property is that one has always
$\dim_{B^{\Gal_E}}D_B(V)\leq \dim_{\Q_p}V$, and the equality holds
iff $V$ is $B$-admissible (See Theorem 2.13 in \cite{FO}). In each
of three cases, one calls $V$ respectively \emph{Hodge-Tate},
\emph{de Rham}, \emph{crystalline} and \emph{semi-stable} (or
\emph{log crystalline}) when $V$ is $B_{HT}$ (resp. $B_{dR}$,
$B_{crys}$, $B_{st}$)-admissible. Furthermore, if the
$B$-admissibility of $V$ holds only for an open subgroup of
$\Gal_E$, then one calls $V$ potentially $B$-admissible.\\

The following simple lemma precedes a further discussion.
\begin{lemma}\label{restriction to an open subgroup}
Let $V$ be a $B$-admissible representation of $\Gal_E$ and
$E'\subset \bar \Q_p$ a finite field extension of $E$. Regarded as
$\Gal_{E'}$-representation via restriction, $V$ is then still
$B$-admissible. Moreover for a crystalline representation $V$, one
has the equality of filtered $\phi$-modules:
$$
(D_{crys,E}(V),\phi,Fil)=(D_{crys,E'}(V)\otimes_{E'_0}E_0,\phi\otimes
\sigma,Fil).
$$
\end{lemma}
\begin{proof}
By $D_{crys,E}=D_{crys,E'}^{\Gal(E'|E)}$ and the injectivity of the
map
$$
\alpha: D_{crys,E'}^{\Gal(E'|E)}\otimes_{E'_0}E_0\to D_{crys,E},
$$
one finds that $V$ is crystalline as $\Gal_{E'}$-representation and
$\alpha$ is an isomorphism.
\end{proof}
Let $r\in \N$. We recall the original definition of
$\Q_{p^r}$-representation in \S7.3 \cite{FO}.
\begin{definition}
A $\Q_{p^r}$-representation of $\Gal_E$ is a finite dimensional
$\Q_{p^r}$-vector space $V$ equipped with a continuous action
$$
\Gal_E\times V\to V
$$
satisfying
$$
g(v_1+v_2)=g(v_1)+g(v_2),\ g(\lambda v)=g(\lambda)g(v)
$$
where $g\in \Gal_E$, and $\lambda\in \Q_{p^r},\ v,v_1,v_2\in V$
\end{definition}
Note after choosing a $\Q_{p^r}$-basis of $V$, one gets a map
$$
f: \Gal_E\to \GL_n(\Q_{p^r}).
$$
However the map $f$ is generally \emph{not} a group homomorphism.
Rather it gives an element in
$$
Z^1_{cts}(\Gal_E,\GL_n(\Q_{p^r}))=\{f(g_1g_2)=f(g_1)g_1(f(g_2)), \
f\ \rm{continuous}\},
$$
where $\Gal_E$ acts on $\GL_n(\Q_{p^r})$ by acting on the entries.
Its isomorphism class is an element in the non-abelian group
cohomology $H^1_{cts}(\Gal_E,\GL_n(\Q_{p^r}))$. Note that the action
of $\Gal_E$ on $\Q_{p^r}$ factors through the quotient $\Gal_{k_E}=
<\phi^f_E>$ where $k_E$ is the residue field of $E$ and
$f_E=[k_E:\F_p]$ is the residue degree of $E$. This implies that the
image of $\Gal_E$ on $\Gal(\Q_{p^r}|\Q_{p})$ is generated by
$\sigma^{f_E\mod r}$. Thus one finds that the action of $\Gal_E$ on
$\Q_{p^r}$ is trivial iff $f_E$ is a multiple of $r$, or
equivalently $E$ contains $\Q_{p^r}$. This fact will simplify our
consideration. The following lemma follows also easily from the
discussion.
\begin{lemma}
Let $(V,\rho)$ be a $\Q_p$-representation of $\Gal_E$. Let $r\in \N$
such that $\Q_{p^r}\subset E$. Then $(V,\rho)$ is a
$\Q_{p^r}$-representation iff there is an injection of
$\Q_p$-algebras $\Q_{p^r}\hookrightarrow \End_{\mathfrak{g}}(V)$,
where $\mathfrak{g}$ is the Lie algebra of the $p$-adic Lie group
$\rho(\Gal_E)\subset \GL(V)$.
\end{lemma}
\begin{definition}
A $\Q_{p^r}$-representation $V$ is called $B$-admissible if it is
$B$-admissible as $\Q_p$-representation.
\end{definition}
For a $\Q_{p^r}$-representation $V$ and each $0\leq m\leq r-1$, one
puts
$$
V_{\sigma^m}:=\Q_{p^r}\otimes_{\sigma^m,\Q_{p^r}}V,
$$
together with tensor product of $\Gal_E$-action. The notation
$\otimes_{\sigma^m}$ signifies the equalities of two tensors:
$$
\lambda(\mu\otimes x)=\lambda\mu\otimes x,\ \lambda\otimes \mu
x=\lambda \mu^{\sigma^m}\otimes x.
$$
\begin{lemma}
Let $f: \Gal_E\to \GL_{n}(\Q_{p^r})$ be an element in
$Z^1_{cts}(\Gal_E,\GL_n(\Q_{p^r}))$ corresponding to $V$. Then the
map $f^{\sigma^m}$, which is defined by
$$
g\mapsto (f(g))^{\sigma^m},
$$
is again an element in $Z^1_{cts}(\Gal_E,\GL_n(\Q_{p^r}))$, and in
fact it corresponds to $V_{\sigma^m}$.
\end{lemma}
\begin{proof}
Let $e_1,\cdots,e_n$ be a basis of $V$. $\tilde e_i:=1\otimes e_i$
basis of $V_{\sigma^m}$. $g(e)=f(g)e$. So
$$
g(\tilde e)=1\otimes g(e)=1\otimes f(g)e=((f(g))^{\sigma^m}1)\otimes
e=f^{\sigma^m}(g)\tilde e.
$$
\end{proof}
\begin{lemma}\label{sigma conjuate is isomorphic to itself}
For each $m$, $V_{\sigma^m}$ is isomorphic to $V$ as
$\Q_p$-representation. In particular, $V$ is $B$-admissible iff
$V_{\sigma^m}$ is $B$-admissible.
\end{lemma}
\begin{proof}
One verifies that the map
$$
V\to V_{\sigma^m}, \ v\mapsto 1\otimes v
$$
is $\Q_p$-linear isomorphism and $\Gal_E$-equivariant.
\end{proof}
However they are not isomorphic as $\Q_{p^r}$-representations as we
will see. For that one studies the invariants
$\{D_{crys,r}^{(m)}(V)\}_{0\leq m\leq r-1}$ introduced for a
$\Q_{p^r}$-representation $V$.
\begin{definition}
Let $V$ be a $\Q_{p^r}$-representation. For each $0\leq m\leq r-1$,
one defines
$$
D_{B,r}^{(m)}(V):=(B\otimes_{\sigma^m,\Q_{p^r}}V)^{\Gal_E}.
$$
\end{definition}
The following trivial lemma underlies the natural direct sum
decomposition for a $\Q_{p^r}$-representation $V$:
\begin{lemma}\label{restriction of scalar of V}
One has a natural isomorphism of $\Gal_{E}$-representations:
$$
\Q_{p^r}\otimes_{\Q_p}V\simeq
\oplus_{m=1}^{r-1}\Q_{p^r}\otimes_{\sigma^m,\Q_{p^r}}V
$$
\end{lemma}
\begin{proof}
One defines the map
$$
a\otimes v\mapsto (av,\cdots,
a^{\sigma^m}v,\cdots,a^{\sigma^{r-1}}v):=f(a\otimes v).
$$
It is injective and $\Q_{p^r}$-linear and hence bijective. One also
checks that $f(ga\otimes gv)=g(f(a\otimes v))$ which follows from
the trivial fact that $\Gal_E$ acts on $\Q_{p^r}$ via a power of
$\sigma$ as discussed above.
\end{proof}
From the lemma one has for example the natural decomposition:
$$
B_{crys}\otimes_{\Q_p}V\simeq
(B_{crys}\otimes_{\Q_{p^r}}\Q_{p^r})\otimes_{\Q_p}V\simeq
B_{crys}\otimes_{\Q_{p^r}}(\Q_{p^r}\otimes_{\Q_p}V)\simeq
\oplus_{m=0}^{r-1}B_{crys}\otimes_{\sigma^m,\Q_{p^r}}V.
$$
This implies that one has the decomposition of $E_0$-vector spaces:
$$
D_{crys}(V)=\oplus_{m=0}^{r-1}D_{crys,r}^{(m)}(V).
$$
Certainly we shall discuss the relation of this decomposition to the
$\sigma$-linear map $\phi$ and the filtration $Fil$ on
$D_{dR}(V)=D_{crys}(V)\otimes_{E_0}E$.
\begin{lemma}\label{phi-module and direct sum}
The map $\phi$ permutes the direct factors $D_{crys,r}^{(m)}(V)$
cyclically. Consequently, one has the decomposition of
$\phi^r$-modules:
$$
(D_{crys,r}(V),\phi^r)=\oplus_{m=0}^{r-1}(D_{crys,r}^{(m)}(V),\phi^r|_{D_{crys,r}^{(m)}(V)}).
$$
Moreover, each $\phi^r$-submodule $(D_{crys,r}^{(m)}(V),\phi^r)$ has
the same Newton slopes.
\end{lemma}
\begin{proof}
Let $d=b\otimes_{\sigma^m}v\in D^{(m)}$, from the formula
$$
\phi(d)=\phi(b)\otimes_{\sigma^{m+1\mod r}}v,
$$
we see that $\phi(d)\in D^{(m+1\mod r)}$. Since $\phi$ is bijective
(so is $\phi^r$ which maps the summand $D^{(m)}$ to itself) and
$\sigma$-linear (which can be viewed as linear map between two
\emph{different} vector spaces), the map
$$
\phi: D_{crys,r}^{(m)}(V)\to D_{crys,r}^{(m+1 \mod r)}(V)
$$
is bijective and $\sigma$-linear isomorphism.\\

If $(f_E,r)=l$, then one can write $af_E=l+br$ with $a,b\in \N$.
Thus $\phi^{af_E}$ is actually linear and induces an isomorphism
$D_{crys,r}^{(m)}(V)\simeq D_{crys,r}^{(m+l \mod r)}(V)$ of $
\phi^r$-modules. In particular, if $f_E$ and $r$ is coprime, then
each direct factor $(D_{crys,r}^{(m)}(V),\phi^r)$ is isomorphic to
each other. Thus we can choose a finite extension $E'$ of $E$ such
that $(f_{E'},r)=1$, and the equality of Newton slopes of each
factor follows from Lemma \ref{restriction to an open subgroup}.
\end{proof}
From the lemma it is also clear that $V$ is crystalline iff
$\dim_{E_0}D_{crys,r}^{(m)}(V)=\dim_{\Q_{p^r}}V$ for either $m$
holds. There is a natural way to construct a $\phi$-module from a
$\phi^r$-module. One starts with a $\phi^r$-module
$(\triangle,\psi)$ over $E_0$ and considers
$$
D(\triangle,\psi)=\Q_p[t]\otimes_{\Q_p[t^r]}\triangle,
$$
where $\triangle$ as $\Q_p[t^r]$-module is given by
$t^r(x):=\psi(x)$.
\begin{lemma}
$D(\triangle,\psi)$ is naturally a $\phi$-module over $E_0$.
Moreover, $D(D_{crys,r}^{(m)}(V),\phi^r)$ is isomorphic to
$(D_{crys}(V),\phi)$.
\end{lemma}
\begin{proof}
We define a $\sigma$-linear map $\phi$ on $D(\triangle,\psi)$ such
that $\phi^r|\triangle=\psi$. For $y=f(t)\otimes \lambda x\in
D(\triangle,\psi)$, one defines $\phi(y)=tf(t)\otimes \lambda^\sigma
x$. The formula
$$
\phi^r(1\otimes x)=t^{r}\otimes x=1\otimes t^r(x)=1\otimes \psi(x)
$$
shows the requirement. To show the isomorphism in the statement, it
suffices to show for $m=0$. As $E_0$-vector space,
$D(D_{crys,r}^{(0)}(V),\phi^r)$ has a natural decomposition
$$
\Q_p[t]\otimes_{\Q_p[t^r]}D_{crys,r}^{(0)}(V)=D_{crys,r}^{(0)}(V)\oplus\cdots\oplus
t^{r-1}D_{crys,r}^{(0)}(V),
$$
so that if $e=\{e_1,\cdots,e_n\}$ is an $E_0$-basis of
$D_{crys,r}^{(0)}(V)$, then $t^me$ is an $E_0$-basis of
$t^mD_{crys,r}^{(0)}(V)$ for $0\leq m\leq r-1$. Note that $\phi^me$
is an $E_0$-basis of $D_{crys,r}^{(m)}(V)$. We defines an
isomorphism of $E_0$-vector spaces $$f:
\Q_p[t]\otimes_{\Q_p[t^r]}D_{crys,r}^{(0)}(V)\to D_{crys}(V)$$ by
sending $t^me_i$ to $\phi^me_i$ for $1\leq i\leq n$. Then one
verifies that $f$ is indeed an isomorphism of $\phi$-modules by
using the definition of the map $\phi$ on
$D(D_{crys,r}^{(0)}(V),\phi^r)$ and the fact that
$\phi^r|_{D_{crys,r}^{(0)}(V)}$ is the original $\phi^r$-module
structure on $D_{crys,r}^{(0)}(V)$.
\end{proof}
This lemma means that, as far as the $\phi$-module structure on
$D_{crys}(V)$ is concerned, one suffices to study the
$\phi^r$-submodule $D_{crys,r}^{(m)}(V)$. Now we consider the
filtration on $D_{dR}(V)$. Again from Lemma \ref{restriction of
scalar of V}, we know that $(D_{dR}(V),Fil)$ is the direct sum of
filtered submodules. Let
$$
D_{dR,r}^{(m)}(V):=(B_{dR}\otimes_{\sigma^m,\Q_{p^r}}V)^{\Gal_E},
$$
and the induced filtration $Fil_m^i:=Fil^i\cap D_{dR,r}^{(m)}(V)$ on
$D_{dR,r}^{(m)}(V)$. Then one has
$$
(D_{dR}(V),Fil)=\oplus_{m=0}^{r-1}(D_{dR,r}^{(m)}(V),Fil_m).
$$
However we must warn that the Hodge slopes of individual factors
(namely the Hodge slopes defined by $Fil_m$) can be
\emph{different}.
\begin{proposition}\label{relation of conjuates on filtered
phi-modules} Let $V$ be a crystalline $\Q_{p^r}$-representation. For
each $0\leq m\leq r-1$, one has
$$
(D^{(m)}_{crys,r}(V),\phi^r,Fil_m)=(D^{(0)}_{crys,r}(V_{\sigma^m}),\phi^r,Fil_0)
$$
as filtered $\phi^r$-modules.
\end{proposition}
\begin{proof}
By Lemma \ref{sigma conjuate is isomorphic to itself}, $D_{crys}(V)$
is isomorphic to $D_{crys}(V_{\sigma^m})$ as filtered
$\phi$-modules. Now since
$$
D_{crys,r}^{(0)}(V_{\sigma^m})=(B_{crys}\otimes_{\Q_{p^r}}(\Q_{p^r}\otimes_{\sigma^m,\Q_{p^r}}V))^{\Gal_E}=
(B_{crys}\otimes_{\sigma^m,\Q_{p^r}}V)^{\Gal_E}=D_{crys,r}^{(m)}(V),
$$
the equality of the statement follows directly from the previous
discussions.
\end{proof}
\begin{lemma}\label{direct factor}
Let $V_1$ and $V_2$ be two $\Q_{p^r}$-representations. One has a
natural decomposition:
$$
V_1\otimes_{\Q_p}V_2\simeq\oplus_{m=0}^{r-1}V_{1}\otimes_{\sigma^m,\Q_{p^r}}V_2.
$$
In particular, $V_1\otimes_{\Q_{p^r}}V_2$ is a direct factor of
$V_1\otimes_{\Q_p}V_2$.
\end{lemma}
\begin{proof}
It follows from Lemma \ref{restriction of scalar of V}.
\end{proof}
\begin{theorem}\label{structure of D_crys}
Let $V_1$ be a crystalline $\Q_{p^r}$-representation and $V$ a
$\Q_p$-representation such that
$$
V\otimes_{\Q_p} \Q_{p^r}\simeq
V_1\otimes_{\Q_{p^r}}V_{1,\sigma}\otimes_{\Q_{p^r}}\cdots\otimes_{\Q_{p^r}}V_{1,\sigma^{r-1}}.
$$
Then one has the equality of filtered $\phi^r$-modules:
$$
(D_{crys}(V),Fil,\phi^r)=
\otimes_{m=0}^{r-1}(D_{crys,r}^{(m)}(V_1),Fil_m,\phi^r).
$$
If $\Q_{p^r}\subset E$, then the above equality is even an equality
of filtered $\phi$-module.
\end{theorem}
\begin{proof}
Put $\tilde V=V\otimes_{\Q_p} \Q_{p^r}$. One has
$$
D_{crys,r}^{(0)}(\tilde
V)=(B_{crys}\otimes_{\Q_{p^r}}(\Q_{p^r}\otimes_{\Q_p}V))^{\Gal_E}=(B_{crys}\otimes_{\Q_{p}}V)^{\Gal_E}=D_{crys}(V).
$$
On the other hand, the tensor decomposition of $\tilde V$ implies
that
$$
D_{crys,r}^{(0)}(\tilde
V)=\otimes_{m=0}^{r-1}D_{crys,r}^{(0)}(V_{1,\sigma^m}),
$$
and $D_{crys,r}^{(0)}(V_{1,\sigma^m})$ is equal to
$D_{crys,r}^{(m)}(V_1)$ by Proposition \ref{relation of conjuates on
filtered phi-modules} as filtered $\phi^r$-modules. \\

Because of Lemma \ref{phi-module and direct sum},
$\otimes_{m=0}^{r-1}D_{crys,r}^{(m)}(V_1)$ is naturally a
$\phi$-module (with the map $\phi$ induced from the $\phi$-structure
on $D_{crys}(V_1)$). Now if $\Q_{p^r}\subset E$, then $\tilde
V=V^{\oplus r}$ as $\Gal_E$-modules. By Lemma \ref{direct factor},
one sees that $V$ is a direct factor of
$$
V_1\otimes_{\Q_{p}}V_{1,\sigma}\otimes_{\Q_{p}}\cdots\otimes_{\Q_{p}}V_{1,\sigma^{r-1}}=V_1^{\otimes
r}.
$$
Therefore, $D_{crys}(V)$ is a direct factor of
$D_{crys}(V_1^{\otimes r})=D_{crys}(V_1)^{\otimes r}$ as filtered
$\phi$-module. In particular, the $\phi$-module structure on
$D_{crys}(V)$ coincides with the restricted one from
$D_{crys}(V)^{\otimes r}$. Write
$D_{crys}(V_1)=\oplus_{m=0}^{r-1}D_{crys,r}^{(m)}(V_1)$. One sees
immediately that $(\otimes_{m=0}^{r-1}D_{crys,r}^{(m)}(V_1),\phi)$
is a direct (filtered) $\phi$-submodule of $D_{crys}(V_1)^{\otimes
r}$. This concludes the proof.
\end{proof}

Now we consider the case that $V$ is polarisable and of Hodge-Tate
weights $0,1$. Here $V$ is polarisable means that there is a perfect
$\Gal_E$-pairing $V\otimes V\to \Q_p(-1)$. This condition implies
that if $\lambda$ is a Newton (resp. Hodge) slopes of $V$, then
$1-\lambda$ is also a Newton (resp. Hodge) slopes of $V$ with the
\emph{same} multiplicity.
\begin{theorem}\label{possible newton slopes}
Let $V$ be a polarisable crystalline representation with Hodge-Tate
weights $0,1$. If there exists a two dimensional crystalline
$\Q_{p^r}$-representation $V_1$ such that
$$
V\otimes_{\Q_p} \Q_{p^r}\simeq
V_1\otimes_{\Q_{p^r}}V_{1,\sigma}\otimes_{\Q_{p^r}}\cdots\otimes_{\Q_{p^r}}V_{1,\sigma^{r-1}}
$$
holds, then it holds that
\begin{itemize}
    \item [(i)] the Hodge slopes of $V_1$ is $\{2r-1\times 0,
1\times 1\}$,
    \item [(ii)] the Newton slopes of $V_1$ is either
$\{2r\times \frac{1}{2r}\}$ or $\{r\times 0,r\times \frac{1}{r}\}$.
\end{itemize}
Consequently, there are only two possible Newton slopes for $V$:
$\{2^r\times \frac{1}{2}\}$ or $\{1\times 0, \cdots, {r\choose
i}\times \frac{i}{r},\cdots, 1\times 1\}$.
\end{theorem}
\begin{proof}
Since the Hodge slopes of $V$ is $\{n\times 0 , n\times 1\}$, by
Theorem \ref{structure of D_crys}, there exists a unique factor
$D^{(i)}_{crys,r}(V_1)$ with two distinct Hodge slopes $\{0,1\}$ and
the other factors have all Hodge slopes zero. WLOG one can assume
that $D^{(0)}_{crys,r}(V_1)$ has Hodge slopes $\{1\times 0,1\times
1\}$ (and any other factor $\{2\times 0\}$). Summing up the Hodge
slopes of all factors, one obtains the Hodge polygon of
$D_{crys}(V_1)$ as claimed. By the admissibility of filtered
$\phi$-module on $D_{crys}(V_1)$, one finds that the Newton slopes
of it must be of form $\{m_1\times 0, m_2\times \lambda\}$ where
$m_1+m_2=2r$ holds, and $\lambda\in \Q$ satisfies $\lambda m_2=1$.
By Lemma \ref{phi-module and direct sum}, one finds that $r|m_i,\
i=1,2$. So $\frac{m_1}{r}+\frac{m_2}{r}=2$ and $m_2\neq 0$. There are two cases to consider:\\

{\bf Case 1}: $m_1=0$. It implies that $m_2=2r$ and
$\lambda=\frac{1}{2r}$.\\

{\bf Case 2:} $m_1\neq 0$. It implies that $m_1=m_2=r$ and
$\lambda=\frac{1}{r}$.
\end{proof}

\section{Newton polygons of abelian varieties of Mumford's type in characteristic $p$}
In this section we apply the previous results to study the Newton
polygons of abelian varieties of Mumford's type. Let $A$ be an
abelian variety of Mumford's type defined over a number field $K$,
and $p$ a rational prime satisfying Assumption \ref{assumption on
p}.
\begin{lemma}[Noot, Proposition 2.1 \cite{Noot1}]
$A$ has potentially good reduction at all places of $K$ over $p$.
\end{lemma}
\begin{proof}
The proof is that of Noot in the case of Mumford's example. Let us
give a sketch only: Let $\mathfrak p$ be a place of $K$ over $p$.
Choose a rational prime $l\neq p$ over which $F$ is inert. By
Deligne's absolute Hodge Proposition 2.9 \cite{De}, one has the
factorization of the Galois representation
$$
\rho_l: \Gal_{K'}\to G(\Q_l)\to \GL(H^1_{et}(\bar A, \Q_l))
$$
for a finite extension $K'$ of $K$. This implies the tensor
decomposition
$$
\rho_{l}\otimes \bar \Q_l=(\otimes_{i=1}^{d}\rho_i)^{\oplus
2^{\epsilon(D)}} .$$ Let $v'$ be a place of $K'$ over $v$ and
$I_{K',v'}$ the inertia group. The monodromy theorem of Grothendieck
implies that $(\rho_l(\gamma)-id)^2=0$ for each $\gamma\in
I_{K',v'}$. The factorization of $\rho_{l}\otimes \bar \Q_l$ implies
that over $\bar \Q_l$ the matrix $\rho_l(\gamma)$ is conjugate to
$(\otimes_{i=1}^{d}\rho_i(\gamma))^{\oplus 2^{\epsilon(D)}}$. Now
that $l$ is inert in $F$, $G^{ss}(\Q_l)$ is simple, and hence
$\{\rho_i(\gamma)\}_{1\leq i\leq d}$ have the same indices. Then
$\rho_l(\gamma)=id,\ \forall \gamma\in I_{K',v'}$, and by Serre-Tate
it has potentially good reduction.
\end{proof}
Let $\mathfrak p_K$ be a prime of $K$ over $p$. By the above lemma
we assume for simplicity that $A$ has the good reduction $A_k$ at
$\mathfrak{p}_K$. Let $K_{\mathfrak{p}_K}$ be the completion of $K$
at $\mathfrak{p}_K$, and $\rho: \Gal_{K_{\mathfrak{p}_K}}\to
\GL(H^1_{et}(\bar A,\Q_p))$ the $p$-adic representation of the local
Galois group. For simplicity we write $E=K_{\mathfrak{p}_K}$ and
$H_{\Q_p}=H^1_{et}(\bar A,\Q_p)$. It is standard that $H_{\Q_p}$ is
known to be a polarisable crystalline representation of Hodge-Tate
weights $\{0,1\}$.
\begin{proposition}\label{first crystalline tensor decomposition}
Over an open subgroup $\Gal_{E'}$ of $\Gal_E$ the representation
$\rho$ admits a natural tensor decomposition:
$$H_{\Q_p}\simeq (V_{\Q_p}\otimes U_{\Q_p})^{\oplus 2^{\epsilon(D)}},$$
where $U_{\Q_p}$ is an unramified representation. Consequently, both
$V_{\Q_p}$ and $U_{\Q_p}$ are crystalline.
\end{proposition}
\begin{proof}
By Proposition \ref{tensor decomposition} one has a tensor
decomposition of $\rho$ over an open subgroup $\Gal_{E'}$:
$$
H_{\Q_p}\simeq (V_{\Q_p}\otimes U_{\Q_p})^{\oplus 2^{\epsilon(D)}}.
$$
Recall that over $\Q_p$, one has
$$
\Cor_{F|\Q}(D)\otimes_{\Q}\Q_p=\prod_{i=1}^{n}\Cor_{F_i|\Q_p}(D\otimes_{F}F_i).
$$
It induces the presentation of $G(\Q_p)$ as a product:
$$
G(\Q_p)=G_1\times G_2,
$$
where $G_1$ corresponds to the $\Gal_{\Q_p}$-orbit containing $\tau$
and $G_2$ the rest $\Gal_{\Q_p}$-orbits (hence $G_2$ may be
trivial). By construction of $V$ and $U$, it is clear that the
projection of $\rho(\Gal_{E'})$ to $G_1$ (resp. $G_2$)-factor acts
on $U_{\Q_p}$ (resp. $V_{\Q_p}$) trivially. Now consider the
Hodge-Tate cocharacter defined by the Hodge-Tate decomposition of
$H_{\Q_p}$:
$$
\mu_{HT}: \G_m(\C_p)\to G(\C_p)=G_1(\C_p)\times G_2(\C_p),
$$
an the Hodge-de Rham cocharacter
$$
\mu_{HdR}: \G_m(\C)\to G(\C)
$$
associated to the Hodge decomposition of $H^1_{B}(A(\C),\Q)$. Let
$C_{HdR}$ (resp. $C_{HT}$) be the $G(\C)$ (resp.
$G(\C_p)$)-conjugacy class of $\mu_{HdR}$ (resp. $\mu_{HT}$). Then
$C_{HdR}$ is defined over the reflex field of $(G,X)$ which is
$\tau(F)\subset \C$. The point is that there is a comparison (see
Page 94 \cite{Noot2} and references therein):
$$
C_{HT}=C_{HdR}\otimes_{F,\tau}\C_p.
$$
Here one uses
$$\tau: F\to \bar \Q\hookrightarrow \bar \Q_p\subset
\hat{\bar{\Q}}_p=\C_p.$$ Hence one deduces that the projection of
$\mu_{HT}$ to the second factor $G_2$ must be trivial since it is
away from the $\tau$-factor. By S. Sen's theorem (the $\Q_p$-Zariski
closure of $\mu_{HT}$ is equal to $\rho(\Gal_{E'^{ur}})$), the
projection to $G_2$ of the inertia group $\subset \Gal_{E'}$ under
$\rho$ is trivial. Thus $U_{\Q_p}$ is an unramified representation
and it is then crystalline (cf. Proposition 7.12 \cite{FO}). As
$V_{\Q_p}\otimes U_{\Q_p}$ is a direct factor of $H_{\Q_p}$, it is
crystalline, and therefore $V_{\Q_p}$, that is a subobject of
$V_{\Q_p}\otimes U_{\Q_p}\otimes U_{\Q_p}^*$, is also crystalline.
\end{proof}

The following result is known among experts. A variant of it was
communicated by L. Berger to the first named author during the
$p$-adic Hodge theory workshop in ICTP, 2009. The first official
proof should appear in the Ph.D thesis of G. Di Matteo (see also
recent preprint \cite{Ma}). Another proof has been communicated to
us by L. Xiao (see \cite{Xiao}).

\begin{theorem}\label{HT implies de-rham }
Let $V$ and $W$ be two $\Q_{p^r}$-representations of $\Gal_E$. If
$V\otimes_{\Q_{p^r}}W$ is de Rham, and one of the tensor factors is
Hodge-Tate, then each tensor factor is de Rham.
\end{theorem}

Applying Theorem \ref{HT implies de-rham } to the tensor factor
$V_{\Q_p}$ in Proposition \ref{first crystalline tensor
decomposition}, one obtains the following
\begin{proposition}\label{second crystalline tensor decomposition}
Making an additional finite field extension $E'\subset E''$ if
necessary, one has a further decomposition of
$\Gal_{E''}$-representation:
$$
V_{\Q_p}\otimes\Q_{p^r}\simeq V_1\otimes_{\Q_{p^r}}
\cdots\otimes_{\Q_{p^r}} V_{r},
$$
where $\Gal_{E''}$ acts on $\Q_{p^r}$ trivially and $V_{i}$ is the
$\sigma^{i-1}$-conjugate of $V_1$. Then each tensor factor $V_i$ is
potentially crystalline.
\end{proposition}
\begin{proof}
Assume $r=2$ for simplicity. The above tensor decomposition implies
the tensor decomposition $\C_p$-representations:
$$
V_{\Q_p}\otimes_{\Q_p}\C_p\simeq
(V_1\otimes_{\Q_{p^2}}\C_p)\otimes_{\C_p}(V_{1,\sigma}\otimes_{\Q_{p^2}}\C_p).
$$
Since $V_{\Q_p}$ is crystalline, it is Hodge-Tate. It implies that
the Sen's operator $\Theta_{V}$ of $V_{\Q_p}$ is diagonalizable
(over $\C_p$). Let $\Theta_{V_1}$ be the Sen's operator of $V_1$. It
can be written naturally as $\Theta_1\oplus\Theta_{1,\sigma}$ where
$\Theta_{1}$ is associated to $V_1\otimes_{\Q_p^2}\C_p$ and
$\Theta_{1,\sigma}$ to $V_{1,\sigma}\otimes_{\Q_p^2}\C_p$. Thus one
has
$$
\Theta_V=\Theta_{1}\otimes id +id\otimes \Theta_{1,\sigma}.
$$
It implies that $\Theta_{1}$ and $\Theta_{1,\sigma}$ are
diagonalizable. Now consider the eigenvalues of them. For that we
use the relation between the Hodge-Tate cocharacter and the
eigenvalues of the Sen's operator: they are related by the maps
$\log$ and $\exp$. Continue the argument about Hodge-Tate
cocharacter in Proposition \ref{first crystalline tensor
decomposition}. So let $\{\tau=\tau_1,\tau_2\}$ be the
$\Gal_{\Q_p}$-orbit of $\tau$. We can assume that in the above
decomposition the $V_1$-factor corresponds to $\tau$. It follows
that the projection of $\mu_{HT}$ to the $V_{1,\sigma}$-factor is
trivial. This implies that the eigenvalues of $\Theta_{1,\sigma}$
are zero. Particularly they are integral. So are those of
$\Theta_1$. Hence $\Theta_{V_1}$ is diagonalizable with integral
eigenvalues. So $V_1$ is Hodge-Tate, and by Theorem \ref{HT implies
de-rham } it is de Rham. By the $p$-adic monodromy theorem,
conjectured by Fontaine and firstly proved by Berger (see
\cite{Be}), it is potentially log crystalline. One shows further
that it is potentially crystalline. Let $N_V$ (resp. $N_{V_1}$) be
the monodromy operator of $V$ (resp. $V_1$). Then as before, one has
the formulas:
$$
N_{V_1}=N_1+N_{1,\sigma}, \ N_V=N_{1}\otimes id+ id\otimes
N_{1,\sigma}.
$$
Since $V$ is crystalline, $N_V=0$. It implies that
$N_1=N_{1,\sigma}=0$. Hence $N_{V_1}=0$ and $V_1$ is potentially
crystalline.
\end{proof}

Now we can prove the main result of the paper.
\begin{theorem}\label{classification and existence of newton polygon}
Notation as above. Put $\mathfrak{p}=\mathfrak{p}_K\cap \sO_F$ and
$r=[F_{\mathfrak{p}}:\Q_p]$. Then the Newton polygon of $A_k$ is
either $\{2^{d+\epsilon(D)}\times \frac{1}{2}\}$ (i.e.
supersingular) or
$$
\{2^{d-r+\epsilon(D)}\times 0, \cdots, 2^{d-r+\epsilon(D)}\cdot{r
\choose i}\times \frac{i}{r}, \cdots, 2^{d-r+\epsilon(D)}\times 1\}.
$$
\end{theorem}
\begin{proof}
The two decomposition results Proposition \ref{first crystalline
tensor decomposition} and \ref{second crystalline tensor
decomposition} show that the condition of Theorem \ref{possible
newton slopes} is satisfied for a finite field extension of $E$.
Note also that the unramified representation $U_{\Q_p}$ contributes
only the multiplicity $2^{d-r}$ to the Newton polygon. The theorem
follows easily from Theorem \ref{possible newton slopes}.
\end{proof}

The concluding paragraph discusses the existence of the possible
Newton polygons. For the example of Mumford, the existence result
was established by Noot (cf. \S3-5 \cite{Noot2}) by studying the
reductions of CM points in a Mumford's family. As a byproduct he has
also obtained a classification of CM points in this case. There are
divided into two cases: Let $F\subset L$ be a maximal subfield of
$D$. Then $L$ can either be written as $E\otimes_{\Q}F$ ($E$ is
necessary an imaginary quadratic extension of $\Q$) or not. To our
purpose one finds the latter case generalizes, and the resulting
generalization gives the necessary existence result. More precisely,
Proposition 5.2 in loc. cit. provides the maximal subfields in $D$
of the second case with the following freedom: Let $\mathfrak{p}$ be
a prime of $F$ over $p$. Then $L$ can be so chosen that $\mathfrak
p$ is split or inert in $L$. Secondly, Lemma 3.5 and Proposition 3.7
in loc. cit. work verbatim for a general $D$: Only one point shall
be taken care of when the rational number $b$ associated to $D$ is
not a square. In this case, one adds the multiplicity two to the
constructions appeared therein. This step gives us an isogeny class
of CM abelian varieties of Mumford's type. Finally the proof of
Proposition 4.4, or rather the method of computing the Newton
polygon for a CM abelian variety modulo a prime works in general.
Thus one can safely conclude the existence result in the general
case.

\end{document}